\documentclass[12pt]{amsproc}
\usepackage{amsmath}
\usepackage{amsfonts}

\theoremstyle{plain}

\newtheorem{corollary}{Corollary}

\newtheorem{definition}{Definition}
\newtheorem{example}{Example}

\newtheorem{notation}{Notation}
\newtheorem{problem}{Problem}
\newtheorem{proposition}{Proposition}
\newtheorem{remark}{Remark}

\newtheorem{theorem}{Theorem}
\numberwithin{equation}{section}
\input{tcilatex}

\begin{document}
\title[Cohomology groups of integral domains and $\phi $-algebras]{%
Cohomology groups of integral domains and $\phi $-algebras}
\author{Mouadh Akriche}
\address{IPEI Monastir\\
05 Avenue Ibn Eljazzar\\
5019 Monastir\\
Tunisia.}
\email{mouadh\_akriche@yahoo.fr}
\author{Mohamed Ali Toumi}
\address{D\'{e}partement de Math\'{e}matiques \\
Facult\'{e} des Sicences de Bizerte\\
7021, Zarzouna, Bizerte\\
Tunisia}
\email{MohamedAli.Toumi@fsb.rnu.tn}

\begin{abstract}
In this paper, we introduce a new homology theory devoted to the study of
linear operators such as local mutipliers and band preserving operators. The
idea is to study the vanishing homology problem. This enables us to
characterize integral domains in which any local multiplier is a multiplier,
which gives a partial answer to a problem posed by R. V. Kadison [J. Algebra
130, No.2, (1990) 494-509]. Finally, we solve the Wickstead problem
[Compositio Math, 35(3) (1977), 225--238] for the class of Archimedean
unital $f$-algebras.
\end{abstract}

\maketitle

\section{Introduction}

The description of algebraic and order properties of an operator is a
primary focus of algebraic analysis. In a unital Archimedean $f$-algebra, we
regard a multiplier as an orthomorphism (that is an order bounded band
preserving operator).\ Then the behavior of algebraic structure reflects
implicit information on the order structure, which is much more accurate
than the one provided by the purely algebraic study.

In [15], the second author et al proved that, for the class of Freudenthal
vector lattices, any band preserving operator in an orthomorphism. For this
purpose, it is very helpful to find invariants such as homological
invariants. One of the purposes in the computation of cohomology groups is
to establish invariants which may be helpful in the classification of the
objects under consideration.

In this paper, we introduce a homology theory for integral domains and
Archimedean unital $f$-algebras which provides information about the
behavior of local multipliers and band preserving operators. The cohomology
theory of associative algebras has been developed by G. Hochschild [7,8,9]
and the 1-, 2-, and 3-dimensional cohomology groups have been interpreted
with reference to classical notions of structure in his papers.

A continuous operator $T$ on algebra $A$ is \textit{local derivation} if for
each $a\in A$ there is a continuous derivation $D_{a}$ on $A$ with $%
D_{a}\left( a\right) =T\left( a\right) .$ This concept was introduced by R.
V. Kadison [10] who showed that if $A$ is a von Neumann algebra, then all
local derivations are in fact derivations. He mentioned also that the set of
all derivations are the 1-cocycles (with respect to the Hochschild homology)
and he sets up a cohomological program for local operators. More precisely,
he invited people to study local higher cohomology (for example, local
2-cocycles with respect to the Hochschild homology).

In Section 2, we study the following problem: Under what conditions on an
integral domain $A,$ any local multiplier is a multiplier? In addition, we
show that the Kadison problem (for local subspace of $2$-cocycles) and the
local multiplier problem are equivalent. Moreover, we are concerned with
connections between structure of an integral domain $A$ and vanishing of its
cohomology groups. To sum up, we prove that on an integral domain $A,$ any
local multiplier is a multiplier if and only if its cohomology groups
vanishes.

The question of whether a band preserving linear operator on Archimedean
vector lattices is automatically order bounded was posed by A. W. Wickstead
[16]. There are several results that guarantee automatic order boundedness
for band preserving operator acting in concrete classes of vector lattices,
see [2,4,5,6,12,13]. The first example of an unbounded band preserving
linear operator was announced by Y. Abramovich, A. I. Veksler and A. V.
Koldunov [2]. Later, they and P. T. N. Mc Polin and A. W. Wickstead [12]
showed that all band preserving operators in a universally complete vector
lattice $A$ are automatically bounded if and only if $A$ is locally
one-dimensional. Hence the Wickstead problem in the class of universally
complete vector lattices was thus reduced to characterization of locally
one-dimensional vector lattices. This characterization was studied in many
works, see [5,6]. There is now a small body of literature devoted to the
study of the Wickstead problem for the class of Archimedean vector lattices.
In fact, S. J. Bernau [4], P. T. N. Mc Polin and A. W. Wickstead [12] and B.
De Pagter [13] proved, by using algebraic and technical tools, that if $T$
is a band preserving linear operator on an Archimedean vector lattice $A$
and if for every positive sequence $\left( x_{n}\right) $ in $A$ which
converges to zero relatively uniformly, $\underset{n}{\inf }\left\{
\left\vert T\left( x_{n}\right) \right\vert \right\} =0$, then $T$ is order
bounded.

In Section 3, we focus our attention on the Wickstead problem on the class
of Archimedean unital $f$-algebras. More precisely, we prove that if $A$ an
Archimedean unital $f$-algebra, then all band preserving operators on $A$
are automatically order bounded if and only if its cohomology groups
vanishes.

In Section 4, we give some examples of integral domains and Archimedean
unital $f$-algebras with vanishing cohomology groups.

We point out that all proofs are purely order theoretical and algebraic in
nature and furthermore do not involve any analytical means. We take it for
granted that the reader is familiar with the notions of vector lattices (or
Riesz spaces) and operators between them. For terminology, notations and
concepts that are not explained in this paper, one can refer to the standard
monographs [1,3,11,14].

\section{The Kadison Problem}

\subsection{Cochains, Coboundary, Cocycles, Cohomology algebra}

Let $A$ be an integral domain and let us denote by $U_{n}\left( A\right) $
the vector space of the $\left( n+1\right) $-linear operators on $A.$ A
linear mapping of $U_{n}\left( A\right) $ into $U_{n+1}\left( A\right) $
analogous to the coboundary operator and leading to the notion of
"cohomology algebra".

We define a "coboundary operator," $d_{n}$, operating on the set of all
cochains as follows:

\begin{definition}
$d_{n}$ maps $U_{n}\left( A\right) $ linearly into $U_{n+1}\left( A\right) $%
. If $f\in U_{1}\left( A\right) ,$ then 
\begin{equation*}
d_{0}\left( f\right) (x_{1},x_{2}))=f\left( x_{1}x_{2}\right) .
\end{equation*}%
If $\Psi \in U_{2}\left( A\right) ,$ then 
\begin{equation*}
d_{1}\left( \Psi \right) (x_{1},x_{2},x_{3}))=\Psi \left(
x_{1}x_{2},x_{3}\right) -\Psi \left( x_{1}x_{3},x_{2}\right) .
\end{equation*}%
If $\Phi \in U_{2n}\left( A\right) $, then 
\begin{eqnarray*}
d_{2n-1}\left( \Phi \right) (x_{1},x_{2},x_{3},...,x_{2n},x_{2n+1})) &=&\Phi
\left( x_{1}x_{2},x_{3},...,x_{2n-1},x_{2n},x_{2n+1})\right) \\
&&-\Phi \left( x_{1}x_{2},x_{3},...,x_{2n-1},x_{2n+1},x_{2n})\right)
\end{eqnarray*}%
for all $n\geq 2.$ If $\Phi \in U_{2n+1}\left( A\right) $, then 
\begin{equation*}
d_{2n}\left( \Phi \right) (x_{1},x_{2},x_{3},...,x_{2n+1},x_{2n+2}))=%
\underset{\sigma \in S_{2n+2}}{\sum }\Phi \left( x_{\sigma \left( 1\right)
}x_{\sigma \left( 2\right) },x_{\sigma \left( 3\right) },...,x_{\sigma
\left( 2n+1\right) },x_{\sigma \left( 2n+2\right) })\right)
\end{equation*}%
for all $n\geq 1.$
\end{definition}

By a simple calculation, we deduce that $d_{n+1}od_{n}=0$, for all $n\geq 0$%
. We now make the customary definitions:

\begin{definition}
A cochain f is called a n-cocycle if $d_{n}\left( f\right) $ $=0$. f is said
to be a n-coboundary if there exists a $\left( n-1\right) $-cochain $g$ such
that $f=d_{n}\left( g\right) $.
\end{definition}

\begin{definition}
For $n\geq 0$, the n-dimensional cohomology group of $A$ denoted by $%
H_{n}(A),$ is the group of the $\left( n+1\right) $-cocycles of $A$ modulo
the subgroup of $n$-coboundaries.
\end{definition}

\subsection{The main results}

A continuous operator $T$ on algebra $A$ is \textit{local derivation} if for
each $a\in A$ there is a continuous derivation $D_{a}$ on $A$ with $%
D_{a}\left( a\right) =T\left( a\right) .$ This concept was introduced by R.
V. Kadison [10] who showed that if $A$ is a von Neumann algebra, than all
local derivations are in fact derivations. He set up a cohomological program
for local operators. More precisely, he invited people to study local higher
cohomology ( for example, local 2-cocycles with respect to the Hochschild
homology). Hence by a simple calculation, for the class of an integral
domain $A,$ the $2$-cocycles, with respect to the Hochschild cohomology, is
the linear space consisting of those bilinear mappings $\Psi $ such that

\begin{equation*}
a\Psi \left( b,c\right) +\Psi \left( a,bc\right) -\Psi \left( ab,c\right)
-c\Psi \left( a,b\right) =0\qquad \left( a,b,c\in A\right) .
\end{equation*}%
By a simple calculation , we deduce that 
\begin{equation*}
\Psi \left( e,a\right) =a\Psi \left( e,a\right) \text{ and }\Psi \left(
a,e\right) =a\Psi \left( e,e\right) \qquad \left( a\in A\right) ,
\end{equation*}%
where $e$ is the unit element of $A.$ Therefore we remark that the set $%
M_{2}\left( A\right) $ of all $\Psi :A\times A\rightarrow A$ such that 
\begin{equation*}
\Psi \left( a,b\right) =b\Psi \left( a,e\right) =a\Psi \left( e,b\right)
=ab\Psi \left( e,e\right) \qquad \left( a,b\in A\right)
\end{equation*}%
is a distinguished sub-space of $2$-cocycles with respect to the Hochschild
homology. Moreover, using the same argument, we also deduce that the set $%
M_{2n}\left( A\right) $ of all $2n$-linear mappings $\Psi :A^{2n}\rightarrow
A$ such that 
\begin{equation*}
\Psi \left( a_{1},a_{2},...,a_{2n}\right) =\left( \underset{1\leq i\leq 2n}{%
\prod }a_{i}\right) \Psi \left( e,e,..,e\right) \qquad \left(
a_{1},a_{2},...,a_{2n}\in A\right)
\end{equation*}%
is a distinguished sub-space of $2n$-cocycles with respect to the Hochschild
homology, for all $n\geq 1.$

\begin{definition}
Let $A$ be an integral domain. A n-linear map $\Psi :A^{n}\rightarrow A$ is
said $n$-multiplier if the following linear mappings 
\begin{equation*}
\Psi ^{i}a_{1},a_{2},..,a_{n-1}:A\rightarrow A;x\mapsto \Psi \left(
a_{1},a_{2},..,a_{i-1},x,a_{i+1}...,a_{n-1}\right) \quad \left( 2\leq i\leq
n-2\right) ,
\end{equation*}%
\begin{equation*}
\Psi ^{1}a_{1},a_{2},..,a_{n-1}:A\rightarrow A;x\mapsto \Psi \left(
x,a_{1},a_{2},..,a_{n-1}\right)
\end{equation*}%
and%
\begin{equation*}
\Psi ^{n}a_{1},a_{2},..,a_{n-1}:A\rightarrow A;x\mapsto \Psi \left(
a_{1},a_{2},..,a_{n-1},x\right)
\end{equation*}%
are multipliers, for all $a_{1},a_{2},..,a_{n-1}\in A$.
\end{definition}

\begin{definition}
Let $A$ be an integral domain. A n-linear map $\Psi :A^{n}\rightarrow A$ is
said local $n$-multiplier if for each $a_{1},a_{2},..,a_{n}\in A$ there
exists a $n$-multiplier linear mapping $\Psi a_{1},a_{2},..,a_{n}$
(depending on $a_{1},a_{2},..,a_{n}$) such that $\Psi
a_{1},a_{2},..,a_{n-1}\left( a_{1},a_{2},..,a_{n}\right) =\Psi \left(
a_{1},a_{2},..,a_{n}\right) .$
\end{definition}

We remark that the up-cited problem posed by R.V. Kadison according to $%
M_{2}\left( A\right) $ is equivalent to the following: Under what conditions
an integral domain $A$ with satisfies the property that any local $2$%
-multiplier on $A^{2}$ is $2$-multiplier? This question motives the
following definition:

\begin{definition}
An algebra $A$ is called a Kadison algebra if any local $2$-multiplier on $%
A^{2}$ is $2$-multiplier.
\end{definition}

The aim of this sub-section is to characterize Kadison integral domains by
using homological approach. In order to hit this mark, we need some
prerequisites.

\begin{notation}
Let $A$ be an integral domain. Let us denote by $C_{n}\left( A\right) $ the
vector space of the $\left( n+1\right) $-linear operators $\Psi $ on $A$
satisfying the following property: $\Psi \left( I_{1},..,I_{n+1}\right)
\subset \underset{1\leq i\leq n+1}{\prod }I_{i,}$ for all $I_{i}$ ideal of $%
A $ and $C_{0}^{m}\left( A\right) $ denotes the vector space of all
multipliers operators $T$ on $A.$
\end{notation}

\begin{remark}
It is not hard to prove that " the coboundary operator" $d_{n}$ satisfies $%
d_{n}\left( C_{n}\left( A\right) \right) \subset C_{n+1}\left( A\right) .$
Therefore by replacing the cochains of the previous homology by $C_{n}\left(
A\right) $ and by keeping the "coboundary operator" $d_{n}$ we have another
homology and its n-dimensional cohomology group of $A$ denoted by $%
H_{n}^{c}(A).$
\end{remark}

\begin{notation}
Let $A$ be an integral domain. Let us denote by $H_{0}^{mc}(A)$ the quotient
group $\ker d_{1}/\left( \func{Im}g\text{ }\left( d_{0/C_{0}^{m}\left(
A\right) }\right) \right) $, where $d_{0/C_{0}^{m}\left( A\right) }$ is the
restriction of $d_{0}$ to $C_{0}^{m}\left( A\right) .$
\end{notation}

\begin{proposition}
Let $A$ be an integral domain. Then \newline
1- The homomorphism $J:H_{0}^{mc}(A)\rightarrow H_{2}(A)$ defined by $%
J\left( \Psi \right) =\Psi _{1}$ with $\Psi _{1}\left(
x_{1},...,x_{4}\right) =\underset{\sigma \in S_{3}}{\sum }x_{1}x_{\sigma
\left( 2\right) }\Psi \left( x_{\sigma \left( 3\right) },x_{\sigma \left(
4\right) }\right) ,$ for any $\Psi \in \overline{\Psi },$ where $\overline{%
\Psi }$ is the equivalent class of $\Psi ,$ is one-to-one.\newline
2- The homomorphism $J^{\prime }:H_{0}^{mc}(A)\rightarrow H_{2}^{c}(A)$
defined by $J^{\prime }\left( \Psi \right) =\Psi _{1}$ with $\Psi _{1}\left(
x_{1},...,x_{4}\right) =\underset{\sigma \in S_{3}}{\sum }x_{1}x_{\sigma
\left( 2\right) }\Psi \left( x_{\sigma \left( 3\right) },x_{\sigma \left(
4\right) }\right) ,$ for any $\Psi \in \overline{\Psi },$ where $\overline{%
\Psi }$ is the equivalent class of $\Psi ,$ is one-to-one.
\end{proposition}

\begin{proof}
1) The mapping $J$ is clearly linear and maps cocycles into cocycles, and
coboundaries into coboundaries. Therefore, it defines a homomorphism of $%
H_{0}^{mc}(A)$ on $H_{2}(A).$ It remains to prove that $J$ is injective. To
this end, let $\Psi \in \overline{\Psi },$ where $\overline{\Psi }\in
H_{0}^{mc}(A)$ such that $\overline{J\left( \Psi \right) }=\overline{0}$ in $%
H_{2}(A).$ Hence $J\left( \Psi \right) =\Psi _{1}$ is symmetric. Then 
\begin{equation}
\Psi _{1}\left( a,b,c,d\right) =\Psi _{1}\left( b,a,c,d\right) \qquad \left(
a,b,c,d\in A\right)  \tag{1}
\end{equation}%
Since $\Psi \in \overline{\Psi },$ it follows that $\Psi $ is symmetric.
Therefore the Equality (1) can be expressed as follows:%
\begin{equation}
a\left[ 2b\Psi \left( c,d\right) +2c\Psi \left( b,d\right) +2d\Psi \left(
b,c\right) \right] =b\left[ 2a\Psi \left( c,d\right) +2c\Psi \left(
a,d\right) +2d\Psi \left( a,c\right) \right]  \tag{2}
\end{equation}%
In Equality (2), if $b=c=d,$ we have%
\begin{equation}
2ab\Psi \left( b,b\right) +2ab\Psi \left( b,b\right) =4b^{2}\Psi \left(
a,b\right)  \tag{3}
\end{equation}%
Since $A$ is an integral domain, we deduce that 
\begin{equation}
a\Psi \left( b,b\right) =b\Psi \left( a,b\right)  \tag{4}
\end{equation}%
Replacing in Equality (4) $b$ by $b+c,$ we get 
\begin{equation*}
a\Psi \left( b+c,b+c\right) =\left( b+c\right) \Psi \left( a,b+c\right) .
\end{equation*}%
Hence 
\begin{equation}
a\Psi \left( b,b\right) +a\Psi \left( c,c\right) +2a\Psi \left( b,c\right)
=b\Psi \left( a,b\right) +b\Psi \left( a,c\right) +c\Psi \left( a,b\right)
+c\Psi \left( a,c\right)  \tag{5}
\end{equation}%
By a simple combination between Equality (4) and Equality (5), we have 
\begin{equation*}
2a\Psi \left( b,c\right) =b\Psi \left( a,c\right) +c\Psi \left( a,b\right)
\end{equation*}%
that is 
\begin{equation}
a\Psi \left( b,c\right) =\frac{1}{2}\left[ b\Psi \left( a,c\right) +c\Psi
\left( a,b\right) \right]  \tag{6}
\end{equation}%
By Equality (6), 
\begin{equation}
c\Psi \left( a,b\right) =\frac{1}{2}\left[ b\Psi \left( a,c\right) +a\Psi
\left( b,c\right) \right] .  \tag{7}
\end{equation}%
By a simple combination between Equality (6) and Equality (7), we have 
\begin{equation}
a\Psi \left( b,c\right) =b\Psi \left( a,c\right)  \tag{8}
\end{equation}%
Then 
\begin{equation*}
\Psi \left( a,b\right) =ab\Psi \left( e,e\right) =a\Psi \left( e,b\right)
=b\Psi \left( a,e\right) ,
\end{equation*}%
where $e$ is the unit element of $A$. Therefore $\Psi $ is $2$-multiplier.
Hence $J$ is one-to-one, which gives the desired result.\newline
2) Using the same argument, we get 2).
\end{proof}

\begin{proposition}
Let $A$ be an integral domain. Then \newline
1- The homomorphism $K:H_{0}^{mc}(A)\rightarrow H_{1}(A)$ defined by $%
K\left( \Psi \right) =\Psi _{1}$ with 
\begin{equation*}
\Psi _{1}\left( x_{1},x_{2},x_{3}\right) =x_{1}\Psi \left(
x_{2},x_{3}\right) -x_{2}\Psi \left( x_{1},x_{3}\right) ,
\end{equation*}%
for any $\Psi \in \overline{\Psi },$ where $\overline{\Psi }$ is the
equivalent class of $\Psi ,$ is one-to-one.\newline
2- The homomorphism $K^{\prime }:H_{0}^{mc}(A)\rightarrow H_{1}^{c}(A)$
defined by $K^{\prime }\left( \Psi \right) =\Psi _{1}$ with 
\begin{equation*}
\Psi _{1}\left( x_{1},x_{2},x_{3}\right) =x_{1}\Psi \left(
x_{2},x_{3}\right) -x_{2}\Psi \left( x_{1},x_{3}\right) ,
\end{equation*}%
for any $\Psi \in \overline{\Psi },$ where $\overline{\Psi }$ is the
equivalent class of $\Psi ,$ is one-to-one.
\end{proposition}

\begin{proof}
1) The mapping $K$ is clearly linear and maps cocycles into cocycles wich is
compatible with coboundaries. Therefore, it defines a homomorphism of $%
H_{0}^{mc}(A)$ on $H_{1}(A).$ It remains to prove that $K$ is injective. To
this end, let $\Psi \in \overline{\Psi },$ where $\overline{\Psi }\in
H_{0}^{mc}(A)$ such that $\overline{K\left( \Psi \right) }=\overline{0}$ in $%
H_{1}(A).$ Since $\Psi \in \overline{\Psi },$ then $\Psi $ is symmetric.
Since $\overline{\Psi }\in H_{0}^{mc}(A),$ then 
\begin{equation*}
\Psi \left( ab,c\right) =\Psi \left( ac,b\right) \qquad \left( a,b,c\in
A\right) .
\end{equation*}%
Hence%
\begin{equation*}
\Psi \left( a,b\right) =\Psi \left( ab,e\right) \qquad \left( a,b\in
A\right) ,
\end{equation*}%
where $e$is the unit element of $A.$ It follows that 
\begin{equation*}
\Psi _{1}\left( a,b,c\right) =a\Psi \left( b,c\right) -b\Psi \left(
a,c\right) =a\Psi \left( bc,e\right) -b\Psi \left( ac,e\right) .
\end{equation*}%
Since $\overline{K\left( \Psi \right) }=\overline{0}$ in $H_{1}(A),$ it
follows that $\Psi _{1}$ satisfies the following property:%
\begin{equation*}
\Psi _{1}\left( a,b,c\right) =-\Psi _{1}\left( a,c,b\right) .
\end{equation*}%
Therefore 
\begin{equation*}
a\Psi \left( bc,e\right) -b\Psi \left( ac,e\right) =-\left( a\Psi \left(
bc,e\right) -c\Psi \left( ab,e\right) \right) .
\end{equation*}%
It follows , if $a=e$ and $b=c$, that 
\begin{equation*}
\Psi \left( b^{2},e\right) -b\Psi \left( b,e\right) =-\left( \Psi \left(
b^{2},e\right) -b\Psi \left( b,e\right) \right) .
\end{equation*}%
Hence 
\begin{equation}
\Psi \left( b,b\right) =\Psi \left( b^{2},e\right) =b\Psi \left( b,e\right) 
\tag{9}
\end{equation}%
Replacing in Equality (9), $a$ by $a+b$, we deduce that 
\begin{equation}
\Psi \left( a,b\right) =\frac{1}{2}\left[ a\Psi \left( e,a\right) +b\Psi
\left( a,e\right) \right]  \tag{10}
\end{equation}%
Replacing in Equality (10) $b$ by $e$, we deduce that 
\begin{equation}
\Psi \left( a,e\right) =a\Psi \left( e,e\right)  \tag{11}
\end{equation}%
Moreover, in Equality (11), if we replace $a$ by $ab$, we have%
\begin{equation*}
\Psi \left( a,b\right) =\Psi \left( ab,e\right) =ab\Psi \left( e,e\right) .
\end{equation*}%
Therefore $\Psi $ is 2-multiplier. Hence $K$ is one-to-one and the proof is
complete.\newline
2) Using the same argument, we get 2).
\end{proof}

Using the same argument we can deduce the following:

\begin{corollary}
Let $A$ be an integral domain. Then \newline
1- The homomorphism $J_{2n}:H_{0}^{mc}(A)\rightarrow H_{2n}(A)$ $(n\geq 1)$
defined by $J_{2n}\left( \Psi \right) =\Psi _{1}$ with 
\begin{equation*}
\Psi _{1}\left( x_{1},...,x_{2n+2}\right) =\underset{\sigma \in S_{3}}{\sum }%
x_{1}x_{\sigma \left( 2\right) }...x_{\sigma \left( 2n-1\right) }\Psi \left(
x_{\sigma \left( 2n+1\right) },x_{\sigma \left( 2n+2\right) }\right) ,
\end{equation*}%
for any $\Psi \in \overline{\Psi },$ where $\overline{\Psi }$ is the
equivalent class of $\Psi ,$ is one-to-one.\newline
2- The homomorphism $J_{2n}^{\prime }:H_{0}^{mc}(A)\rightarrow H_{2n}^{c}(A)$
defined by $J_{2n}^{\prime }\left( \Psi \right) =\Psi _{1}$ with 
\begin{equation*}
\Psi _{1}\left( x_{1},...,x_{2n+2}\right) =\underset{\sigma \in S_{3}}{\sum }%
x_{1}x_{\sigma \left( 2\right) }...x_{\sigma \left( 2n-1\right) }\Psi \left(
x_{\sigma \left( 2n+1\right) },x_{\sigma \left( 2n+2\right) }\right) ,
\end{equation*}%
for any $\Psi \in \overline{\Psi },$ where $\overline{\Psi }$ is the
equivalent class of $\Psi ,$ is one-to-one.\newline
3- The homomorphism $J_{2n-1}:H_{0}^{mc}(A)\rightarrow H_{2n-1}(A)$ defined
by $J_{2n-1}\left( \Psi \right) =\Psi _{1}$ with 
\begin{equation*}
\Psi _{1}\left( x_{1},...,x_{2n+1}\right) =\underset{1\leq i\leq 2n-2}{\prod
x_{i}}\left( x_{2n-1}\Psi \left( x_{2n},x_{2n+1}\right) -x_{2n}\Psi \left(
x_{2n-1},x_{2n+1}\right) \right) ,
\end{equation*}%
for any $\Psi \in \overline{\Psi },$ where $\overline{\Psi }$ is the
equivalent class of $\Psi ,$ is one-to-one.\newline
4- The homomorphism $J_{2n-1}^{\prime }:H_{0}^{mc}(A)\rightarrow
H_{2n-1}^{c}(A)$ defined by $J_{2n-1}^{\prime }\left( \Psi \right) =\Psi
_{1} $ with 
\begin{equation*}
\Psi _{1}\left( x_{1},...,x_{2n+1}\right) =\underset{1\leq i\leq 2n-2}{\prod
x_{i}}\left( x_{2n-1}\Psi \left( x_{2n},x_{2n+1}\right) -x_{2n}\Psi \left(
x_{2n-1},x_{2n+1}\right) \right) ,
\end{equation*}%
for any $\Psi \in \overline{\Psi },$ where $\overline{\Psi }$ is the
equivalent class of $\Psi ,$ is one-to-one.
\end{corollary}

We have gathered all ingredients for the main results of this section.

\begin{theorem}
Let $A$ be an integral domain. Then the following assertions are equivalent.%
\newline
i) $H_{0}^{mc}(A)=\left\{ 0\right\} $\newline
ii) $H_{n}^{c}(A)=\left\{ 0\right\} $, for all $n\geq 1$\newline
iii) There exists $n\geq 1$ such that $H_{n}^{c}(A)=\left\{ 0\right\} $%
\newline
vi) $A$ is a Kadison space.
\end{theorem}

\begin{proof}
$i)\Rightarrow ii)$ Let $T$ be a local multiplier on $A.$ Let us define the
following $2$-local multiplier $\Psi $ on $A$ by $\Psi \left( a,b\right)
=T\left( ab\right) .$ Since $H_{0}^{mc}(A)=\left\{ 0\right\} $ then $T$ is a
multiplier. Therefore $T\left( a\right) =aT\left( e\right) .$ Using the same
argument, we deduce that 
\begin{equation*}
\Psi \left( x_{1},..,x_{n}\right) =\left( \underset{1\leq i\leq n}{\prod }%
x_{i}\right) \Psi \left( e,..,e\right) \qquad \left( x_{1},..,x_{n}\in
A\right)
\end{equation*}%
for all $\Psi \in C_{n}\left( A\right) $, where $e.$is the unit element of $%
A $. It follows that $H_{n}^{c}(A)=\left\{ 0\right\} $, for all $n\geq 1.$

$ii)\Rightarrow iii)$ Obvious.

$iii)\Rightarrow vi)$ If there exists $n\geq 1$ such that $%
H_{n}^{c}(A)=\left\{ 0\right\} ,$ then by the previous propositions $%
H_{0}^{mc}(A)=\left\{ 0\right\} $ then any local multiplier on $A$ is a
multiplier. Therefore $A$ is a Kadison space.

$vi\Rightarrow i)$ $A$ is a Kadison space, then any local multiplier on $A$
is a multiplier. Let $\Psi $ be $2$-local multiplier $\Psi $ on $A.$ Then $%
\Psi \left( a,.\right) :A\rightarrow A$ defined by $\Psi \left( a,.\right)
\left( b\right) =\Psi \left( a,b\right) $ is a local multiplier. By the fact
that $A$ is a Kadison space then $\Psi \left( a,b\right) =b\Psi \left(
a,e\right) ,$ for all $ab,\in A.$ Using the same idea, we deduce that $\Psi
\left( a,b\right) =ab\Psi \left( e,e\right) =a\Psi \left( e,b\right) =b\Psi
\left( a,e\right) $ and we are done.
\end{proof}

\begin{corollary}
Let $A$ be an integral domain. If there exists $n\geq 1$ such that $%
H_{n}(A)=\left\{ 0\right\} $ then $A$ is a Kadison space.
\end{corollary}

We end this section with the following problem:

\begin{problem}
Is it true, for the case of an integral domain, that if $A$ is a Kadison
space then there exists $n\geq 1$ such that $H_{n}(A)=\left\{ 0\right\} ?$
\end{problem}

\section{The Wickstead Problem}

\subsection{Definitions and notations}

In order to avoid unnecessary repetition we will suppose that all vector
lattices and $\ell $-algebras under consideration are \textbf{Archimedean}.

Let us recall some of the relevant notions. Let $A$ be a vector lattice. A
linear operator $T:A\rightarrow A$ is called \textit{band preserving } if $%
T\left( x\right) \perp y$ whenever $x\perp y$ in $A$. A linear mapping 
\textit{T} $\in \mathcal{L}(A,B)$ is called \textit{order bounded } if $T$
maps order bounded subsets of $A$ onto order bounded subsets of $B.$ An
order bounded band preserving on\textit{\ }$A$ is called \textit{%
orthomorphism}.

A bilinear operator $\Psi :A\times A\rightarrow A$ is called \textit{%
separately band preserving} (resp \textit{separately ideal preserving})
provided that the following mappings%
\begin{equation*}
\Psi (.,x):y\mapsto \Psi (y,x)\text{ and }\Psi (x,.):y\mapsto \Psi (x,y)%
\text{ }\qquad (y\in A)
\end{equation*}%
are band preserving (resp ideal preserving) for all $x\in A$.

In the following lines, we recall definitions and some basic facts about $f$%
-algebras. For more information about this field, one can refer to [1,11]. A
(real) algebra \textit{A} which is simultaneously a vector lattice such that
the partial ordering and the multiplication in $A$ are compatible, that is $%
a,b\in A^{+}$ implies $ab\in A^{+}$ is called \textit{a} \textit{%
lattice-ordered algebra}( briefly $\ell $\textit{-algebra}). In an $\ell $%
-algebra \textit{A }we denote the collection of all nilpotent elements of 
\textit{A} by \textit{N(A)}. An $\ell $-algebra \textit{A }is said to be 
\textit{semiprime} if $N(A)=\{0\}$\textit{. }An $\ell $-algebra \textit{A }%
is called an \textit{f-algebra }if \textit{A} verifies the property that $%
a\wedge b=0$ and $c\geq 0$ imply $ac\wedge b=ca\wedge b=0.$ A unital \textit{%
f-}algebra (i.e., an \textit{f-}algebra with a unit element) is called $\phi 
$\textit{-algebra}.

The vector lattice\textit{\ }$A$ is called \textit{Dedekind }$\sigma $-%
\textit{complete} if for each non-void countable majorized set $B\subset A$, 
$\sup $ $B$ exists in $A$. The vector lattice\textit{\ }$A$ is called 
\textit{laterally complete} provided that every orthogonal system in \textit{%
A} has a supremum in \textit{A.} If $A$ is Dedekind complete and laterally
complete, then $A$ is said to be\textit{\ universally complete}. Every
vector lattice \textit{A} has a \textit{universal completion }$A^{u}$, this
means that there exists a unique (up to a lattice isomorphism) universally
complete vector lattice $A^{u}$ such that \textit{A} can be identified \
with an order dense sublattice of $A^{u}.$ For more properties about
universal completion, see [11, chap VII, section 51].

\subsection{The main results}

The present section considers band preserving linear operators on $\phi $%
\textit{-algebras}. More precisely, we are mainly concerned with
characterizing $\phi $\textit{-algebras} on which any band preserving linear
operator is order bounded.

This question motives the following definition:

\begin{definition}
An $\ell -$algebra $A$ is called a Wickstead algebra if any band preserving
linear operator on $A$ is order bounded.
\end{definition}

The aim of this sub-section is to characterize the Wickstead $\phi $%
-algebras by using homological approach. In order to hit this mark, we need
some prerequisites.

\begin{notation}
Let $A$ be a $\phi $\textit{-algebra}. Let us denote by $F_{n}\left(
A\right) $ the vector space of the $\left( n+1\right) $-linear \textit{%
separately }band preserving operators $\Psi $ on $A$ and $Orth\left(
A\right) $ denotes the vector space of all orthomorphisms $T$ on $A$
\end{notation}

\begin{remark}
It is not hard to prove that "the coboundary operator" $d_{n}$ satisfies $%
d_{n}\left( F_{n}\left( A\right) \right) \subset F_{n+1}\left( A\right) .$
Therefore by replacing the cochains of the previous homology by $F_{n}\left(
A\right) $ and by keeping the "coboundary operator" $d_{n}$ we have another
homology and its n-dimensional cohomology group of $A$ denoted by $%
H_{n}^{o}(A).$
\end{remark}

\begin{notation}
Let $A$ be a $\phi $\textit{-algebra}. Let us denote by $H_{0}^{oo}(A)$ the
quotient group $\ker d_{1}/\left( \func{Im}g\text{ }\left( d_{0/Orth\left(
A\right) }\right) \right) ,$ where $d_{0/Orth\left( A\right) }$ is the
restriction of $d_{0}$ to $Orth\left( A\right) .$
\end{notation}

\begin{proposition}
Let $A$ be a $\phi $\textit{-algebra}. Then \newline
1- The homomorphism $J:H_{0}^{oo}(A)\rightarrow H_{2}(A)$ defined by $%
J\left( \Psi \right) =\Psi _{1}$ with $\Psi _{1}\left(
x_{1},...,x_{4}\right) =\underset{\sigma \in S_{3}}{\sum }x_{1}x_{\sigma
\left( 2\right) }\Psi \left( x_{\sigma \left( 3\right) },x_{\sigma \left(
4\right) }\right) ,$ for any $\Psi \in \overline{\Psi },$ where $\overline{%
\Psi }$ is the equivalent class of $\Psi ,$ is one-to-one.\newline
2- The homomorphism $J^{\prime }:H_{0}^{oo}(A)\rightarrow H_{2}^{o}(A)$
defined by $J^{\prime }\left( \Psi \right) =\Psi _{1}$ with $\Psi _{1}\left(
x_{1},...,x_{4}\right) =\underset{\sigma \in S_{3}}{\sum }x_{1}x_{\sigma
\left( 2\right) }\Psi \left( x_{\sigma \left( 3\right) },x_{\sigma \left(
4\right) }\right) ,$ for any $\Psi \in \overline{\Psi },$ where $\overline{%
\Psi }$ is the equivalent class of $\Psi ,$ is one-to-one.
\end{proposition}

\begin{proof}
1) The mapping $J$ is clearly linear and maps cocycles into cocycles, and
coboundaries into coboundaries. Therefore, it defines a homo- morphism of $%
H_{0}^{oo}(A)$ on $H_{2}(A).$ It remains to prove that $J$ is injective. To
this end, let $\Psi \in \overline{\Psi },$ where $\overline{\Psi }\in
H_{0}^{oo}(A)$ such that $\overline{J\left( \Psi \right) }=\overline{0}$ in $%
H_{2}(A).$ Hence $J\left( \Psi \right) =\Psi _{1}$ is symmetric. Then 
\begin{equation}
\Psi _{1}\left( a,b,c,d\right) =\Psi _{1}\left( b,a,c,d\right) \qquad \left(
a,b,c,d\in A\right)  \tag{12}
\end{equation}%
Since $\Psi \in \overline{\Psi },$ it follows that $\Psi $ is symmetric.
Therefore the Equality (12) can be expressed as follows:%
\begin{equation}
a\left[ 2b\Psi \left( c,d\right) +2c\Psi \left( b,d\right) +2d\Psi \left(
b,c\right) \right] =b\left[ 2a\Psi \left( c,d\right) +2c\Psi \left(
a,d\right) +2d\Psi \left( a,c\right) \right]  \tag{13}
\end{equation}%
In Equality (13), if $b=c=d,$ we have%
\begin{equation}
2ab\Psi \left( b,b\right) +2ab\Psi \left( b,b\right) =4b^{2}\Psi \left(
a,b\right)  \tag{14}
\end{equation}%
Since $A$ is a $\phi $\textit{-algebra} and $\Psi $ is separately band
preserving, we deduce that 
\begin{equation}
a\Psi \left( b,b\right) =b\Psi \left( a,b\right)  \tag{15}
\end{equation}%
Replacing in Equality (14) $b$ by $b+c,$ we get 
\begin{equation*}
a\Psi \left( b+c,b+c\right) =\left( b+c\right) \Psi \left( a,b+c\right) .
\end{equation*}%
Hence 
\begin{equation}
a\Psi \left( b,b\right) +a\Psi \left( c,c\right) +2a\Psi \left( b,c\right)
=b\Psi \left( a,b\right) +b\Psi \left( a,c\right) +c\Psi \left( a,b\right)
+c\Psi \left( a,c\right)  \tag{16}
\end{equation}%
By a simple combination between Equality (15) and Equality (16), we have 
\begin{equation*}
2a\Psi \left( b,c\right) =b\Psi \left( a,c\right) +c\Psi \left( a,b\right)
\end{equation*}%
that is 
\begin{equation}
a\Psi \left( b,c\right) =\frac{1}{2}\left[ b\Psi \left( a,c\right) +c\Psi
\left( a,b\right) \right]  \tag{17}
\end{equation}%
By Equality (6), 
\begin{equation}
c\Psi \left( a,b\right) =\frac{1}{2}\left[ b\Psi \left( a,c\right) +a\Psi
\left( b,c\right) \right] .  \tag{18}
\end{equation}%
By a simple combination between Equality (17) and Equality (18), we have 
\begin{equation}
a\Psi \left( b,c\right) =b\Psi \left( a,c\right)  \tag{19}
\end{equation}%
Then 
\begin{equation*}
\Psi \left( a,b\right) =ab\Psi \left( e,e\right) =a\Psi \left( e,b\right)
=b\Psi \left( a,e\right) ,
\end{equation*}%
where $e$ in the unit element of $A$. Therefore $\Psi $ is 2-multiplier.
Hence $J$ is one-to-one, which completes the proof.\newline
2) Using the same argument, we get 2).
\end{proof}

\begin{proposition}
Let $A$ be a $\phi $\textit{-algebra}. Then \newline
1- The homomorphism $K:H_{0}^{oo}(A)\rightarrow H_{1}(A)$ defined by $%
K\left( \Psi \right) =\Psi _{1}$ with 
\begin{equation*}
\Psi _{1}\left( x_{1},x_{2},x_{3}\right) =x_{1}\Psi \left(
x_{2},x_{3}\right) -x_{2}\Psi \left( x_{1},x_{3}\right) ,
\end{equation*}%
for any $\Psi \in \overline{\Psi },$ where $\overline{\Psi }$ is the
equivalent class of $\Psi ,$ is one-to-one.\newline
2- The homomorphism $K^{\prime }:H_{0}^{oo}(A)\rightarrow H_{1}^{o}(A)$
defined by $K^{\prime }\left( \Psi \right) =\Psi _{1}$ with 
\begin{equation*}
\Psi _{1}\left( x_{1},x_{2},x_{3}\right) =x_{1}\Psi \left(
x_{2},x_{3}\right) -x_{2}\Psi \left( x_{1},x_{3}\right) ,
\end{equation*}%
for any $\Psi \in \overline{\Psi },$ where $\overline{\Psi }$ is the
equivalent class of $\Psi ,$ is one-to-one.
\end{proposition}

\begin{proof}
2) The mapping $K$ is clearly linear and maps cocycles into cocycles, and
coboundaries into coboundaries. Therefore, it defines a homo- morphism of $%
H_{0}^{oo}(A)$ on $H_{1}(A).$ It remains to prove that $K$ is injective. To
this end, let $\Psi \in \overline{\Psi },$ where $\overline{\Psi }\in
H_{0}^{oo}(A)$ such that $\overline{K\left( \Psi \right) }=\overline{0}$ in $%
H_{1}(A).$ Since $\Psi \in \overline{\Psi },$ then $\Psi $ is symmetric.
Since $\overline{\Psi }\in H_{0}^{oo}(A),$ then 
\begin{equation*}
\Psi \left( ab,c\right) =\Psi \left( ac,b\right) \qquad \left( a,b,c\in
A\right) .
\end{equation*}%
It follows that 
\begin{equation*}
\Psi \left( a,b\right) =\Psi \left( ab,e\right) \qquad \left( a,b\in A\right)
\end{equation*}%
where $e$ is the unit element of $A.$Hence 
\begin{equation*}
\Psi _{1}\left( a,b,c\right) =a\Psi \left( b,c\right) -b\Psi \left(
a,c\right) =a\Psi \left( bc,e\right) -b\Psi \left( ac,e\right) .
\end{equation*}%
Since $\overline{K\left( \Psi \right) }=\overline{0}$ in $H_{1}(A),$ it
follows that $K\left( \Psi \right) =\Psi _{1}$ satisfies the following
property: 
\begin{equation*}
\Psi _{1}\left( a,b,c\right) =-\Psi _{1}\left( a,c,b\right) .
\end{equation*}%
Therefore 
\begin{equation*}
a\Psi \left( bc,e\right) -b\Psi \left( ac,e\right) =-\left( a\Psi \left(
bc,e\right) -c\Psi \left( ab,e\right) \right) .
\end{equation*}%
It follows , if $a=e$ and $b=c$, that 
\begin{equation*}
\Psi \left( b^{2},e\right) -b\Psi \left( b,e\right) =-\left( \Psi \left(
b^{2},e\right) -b\Psi \left( b,e\right) \right) .
\end{equation*}%
Hence 
\begin{equation}
\Psi \left( b,b\right) =\Psi \left( b^{2},e\right) =b\Psi \left( b,e\right) 
\tag{20}
\end{equation}%
Replacing in Equality (20), $a$ by $a+b$, we deduce that 
\begin{equation}
\Psi \left( a,b\right) =\frac{1}{2}\left[ a\Psi \left( e,a\right) +b\Psi
\left( a,e\right) \right]  \tag{21}
\end{equation}%
Replacing in Equality (21) $b$ by $e$, we deduce that 
\begin{equation}
\Psi \left( a,e\right) =a\Psi \left( e,e\right)  \tag{22}
\end{equation}%
Moreover, in Equality (22), if we replace $a$ by $ab$, we have%
\begin{equation*}
\Psi \left( a,b\right) =\Psi \left( ab,e\right) =ab\Psi \left( e,e\right) .
\end{equation*}%
Therefore $\Psi $ is 2-multiplier. Hence $K$ is one-to-one and we are done.%
\newline
2) Using the same argument, we get 2).
\end{proof}

Using the same argument we can deduce the following:

\begin{corollary}
Let $A$ be a $\phi $\textit{-algebra}. Then \newline
1- The homomorphism $J_{2n}:H_{0}^{oo}(A)\rightarrow H_{2n}(A)$ $(n\geq 1)$
defined by $J_{2n}\left( \Psi \right) =\Psi _{1}$ with 
\begin{equation*}
\Psi _{1}\left( x_{1},...,x_{2n+2}\right) =\underset{\sigma \in S_{3}}{\sum }%
x_{1}x_{\sigma \left( 2\right) }...x_{\sigma \left( 2n-1\right) }\Psi \left(
x_{\sigma \left( 2n+1\right) },x_{\sigma \left( 2n+2\right) }\right) ,
\end{equation*}%
for any $\Psi \in \overline{\Psi },$ where $\overline{\Psi }$ is the
equivalent class of $\Psi ,$ is one-to-one.\newline
2- The homomorphism $J_{2n}^{\prime }:H_{0}^{oo}(A)\rightarrow H_{2n}^{o}(A)$
defined by $J_{2n}^{\prime }\left( \Psi \right) =\Psi _{1}$ with 
\begin{equation*}
\Psi _{1}\left( x_{1},...,x_{2n+2}\right) =\underset{\sigma \in S_{3}}{\sum }%
x_{1}x_{\sigma \left( 2\right) }...x_{\sigma \left( 2n-1\right) }\Psi \left(
x_{\sigma \left( 2n+1\right) },x_{\sigma \left( 2n+2\right) }\right) ,
\end{equation*}%
for any $\Psi \in \overline{\Psi },$ where $\overline{\Psi }$ is the
equivalent class of $\Psi ,$ is one-to-one.\newline
3- The homomorphism $J_{2n-1}:H_{0}^{oo}(A)\rightarrow H_{2n-1}(A)$ defined
by $J_{2n-1}\left( \Psi \right) =\Psi _{1}$ with 
\begin{equation*}
\Psi _{1}\left( x_{1},...,x_{2n+1}\right) =\left( \underset{1\leq i\leq 2n-2}%
{\prod x_{i}}\right) \left( x_{2n-1}\Psi \left( x_{2n},x_{2n+1}\right)
-x_{2n}\Psi \left( x_{2n-1},x_{2n+1}\right) \right) ,
\end{equation*}%
for any $\Psi \in \overline{\Psi },$ where $\overline{\Psi }$ is the
equivalent class of $\Psi ,$ is one-to-one.\newline
4- The homomorphism $J_{2n-1}^{\prime }:H_{0}^{oo}(A)\rightarrow
H_{2n-1}^{o}(A)$ defined by $J_{2n-1}\left( \Psi \right) =\Psi _{1}$ with 
\begin{equation*}
\Psi _{1}\left( x_{1},...,x_{2n+1}\right) =\left( \underset{1\leq i\leq 2n-2}%
{\prod x_{i}}\right) \left( x_{2n-1}\Psi \left( x_{2n},x_{2n+1}\right)
-x_{2n}\Psi \left( x_{2n-1},x_{2n+1}\right) \right) ,
\end{equation*}%
for any $\Psi \in \overline{\Psi },$ where $\overline{\Psi }$ is the
equivalent class of $\Psi ,$ is one-to-one.
\end{corollary}

We have gathered all ingredients for the main results of this section. The
proof is omitted since it is similar to the proof of the previous theorem.

\begin{theorem}
Let $A$ be a $\phi $\textit{-algebra}. Then the following assertions are
equivalent.\newline
i) $H_{0}^{oo}(A)=\left\{ 0\right\} $\newline
ii) $H_{n}^{o}(A)=\left\{ 0\right\} $, for all $n\geq 1$\newline
iii) There exists $n\geq 1$ such that $H_{n}^{o}(A)=\left\{ 0\right\} $%
\newline
vi) $A$ is a Wickstead space.
\end{theorem}

\begin{corollary}
Let $A$ be a $\phi $\textit{-algebra}. If there exists $n\geq 1$ such that $%
H_{n}(A)=\left\{ 0\right\} $ then $A$ is a Wickstead space.
\end{corollary}

Since in any $\phi $\textit{-algebra} the two notions orthomorphisms and and
multipliers are the same, we deduce the following result.

\begin{corollary}
Let $A$ be a $\phi $\textit{-algebra}. Then a Wickstead space is a Kadison
space.
\end{corollary}

We end this section with the following problem:

\begin{problem}
Is it true, for the case of $\phi $\textit{-algebras}, that if $A$ is a
Wickstead space then there exists $n\geq 1$ such that $H_{n}(A)=\left\{
0\right\} ?$
\end{problem}

\section{Examples}

Recall that a norm $\left\Vert .\right\Vert $ on a vector lattice is said to
be a \textit{lattice norm} whenever $\left\vert x\right\vert \leq \left\vert
y\right\vert $ in $A$ implies $\left\Vert x\right\Vert \leq \left\Vert
y\right\Vert $. A vector lattice equipped with a lattice norm is known as a 
\textit{normed vector lattice.}

It is shown in [1, Theorem 4.76], that any band preserving operator on a
Banach lattice is inevitably order bounded. Hence we deduce that:

\begin{theorem}
Any Banach vector lattice is a Wickstead space.
\end{theorem}

Let $A$ be a vector lattice and let $0\leq a\in A.$ An element $0\leq e\in A$
is called a \textit{component }of $a$ if $e\wedge \left( a-e\right) =0.$

\begin{definition}
A vector lattice $A$ is called a Freudenthal vector lattice (or
Hyper-Archimedean) if $A$ satisfies the following property:\newline
If $0\leq x\leq e$ holds in $A$, then there exist positive real numbers $%
\alpha _{1},....,\alpha _{n}$ and components $e_{1},..,e_{n}$ of $e$
satisfying $x=\underset{1\leq i\leq n}{\sum }\alpha _{i}e_{i}.$
\end{definition}

\begin{example}
The vector space of all real stationary sequences is an atomic Freudenthal
vector lattice.
\end{example}

It is shown in [15], that any band preserving operator on a Freudenthal
vector lattice is inevitably order bounded. Hence we deduce that:

\begin{theorem}
Any Freudenthal vector lattice is a Wickstead space.
\end{theorem}

P. T. N. Mc Polin and A. W. Wickstead [12] showed that all band preserving
operators in a universally complete vector lattice $A$ are automatically
bounded if and only if $A$ is locally one-dimensional. Hence we have:

\begin{theorem}
Any locally one-dimensional vector lattice is a Wickstead space.
\end{theorem}

\end{document}